\documentclass[12pt]{amsart}            

\usepackage{latexsym,amssymb}
\usepackage{amsmath,amsthm}      
\usepackage{amsopn}              
\usepackage{cite}                

\title{Star-operations induced by overrings.}

\author{Sharon M. Clarke}

\newcommand{\nnn}[1]{\eqref{#1}}

\newtheorem{theorem}{Theorem}[section]
\newtheorem{lemma}[theorem]{Lemma}
\newtheorem{proposition}[theorem]{Proposition}
\newtheorem{corollary}[theorem]{Corollary}
\newtheorem{definition}[theorem]{Definition}
\newtheorem{claim}{Claim}

\theoremstyle{definition}
\newtheorem{example}[theorem]{Example}

\newcommand{\M}{\mathcal}
\newcommand{\ra}{\Rightarrow}
\newcommand{\la}{\Leftarrow}

\renewcommand{\a}{\alpha}

\def\<#1,#2>{\langle{#1},{#2}\rangle} 

\def\sideremark#1{\ifvmode\leavevmode\fi\vadjust{\vbox to0pt{\vss
\hbox to 0pt{\hskip\hsize\hskip1em
\vbox{\hsize2cm\tiny\raggedright\pretolerance10000
\noindent #1\hfill}\hss}\vbox to8pt{\vfil}\vss}}}

\begin{document}

\begin{abstract}

Let $D$ be an integral domain with quotient field $K$. A star-operation $\star$ on $D$ is 
a closure operation $A \longmapsto 
A^\star$ on the set of nonzero fractional ideals, $F(D)$, of $D$ satisfying the 
properties: $(xD)^\star = xD$ and $(xA)^\star = xA^\star$ for all $x \in K^\ast$ and 
$A \in F(D)$. Let ${\M S}$ be a multiplicatively closed set of ideals
of $D$. For $A \in F(D)$ define $A_{\M S} = \{x \in K \mid xI \subseteq{A}$, for some $I 
\in {\M S}\}$. Then $D_{\M S}$ is an overring of $D$ and
$A_{\M S}$ is a fractional ideal of $D_{\M S}$. Let ${\M S}$ be a multiplicative set of 
finitely generated nonzero ideals of $D$ and $A \in F(D)$, then the map 
$A \longmapsto A_{\M S}$ is a finite character star-operation if and only if for each 
$I \in {\M S}$, $I_v = D$. We give an example to show that this result is not true if the 
ideals are not assumed to be finitely generated. In general, the map $A \longmapsto 
A_{\M S}$ is a star-operation if and only if $\bar {\M S}$, the saturation of ${\M S}$, is
a localizing GV-system. We also discuss star-operations given of the form $A \longmapsto
\cap AD_\alpha$, where $D = \cap D_\alpha$.
  
\end{abstract}

\maketitle  


\pagestyle{myheadings}
\markboth{Sharon \ M. \ Clarke}{Star-operations induced by overrings}


\section{Introduction}

Let $D$ be an integral domain with quotient field $K$. Let $K^\ast = K - \{0\}$.
A {\it fractional ideal} 
$A$ of $D$ is a $D$-submodule of $K$ such that $xA \subseteq{D}$ for some 
nonzero $x \in D$. Note that the usual ideals of $D$, or sometimes called 
{\it integral ideals}, are
also fractional ideals. Let $F(D)$ be the set of nonzero fractional ideals of 
$D$ and $f(D)$ be the subset of $F(D)$ consisting of the nonzero finitely 
generated fractional ideals of $D$. Let $I(D)$ denote the set of nonzero 
integral ideals of $D$. 

\begin{definition}\label{star}

A {\bf star-operation} on the set F(D) is a mapping $\star : F(D) 
\longrightarrow F(D)$ such that for all $A, B \in F(D)$, and for all $a \in
K^\ast$:

(1) $(a)^\star = (a)$ and $(aA)^\star = aA^\star$,

(2) $A \subseteq{A^\star}$, and $A \subseteq{B}$ implies that $A^\star 
    \subseteq{B^\star}$, and
    
(3) $(A^\star)^\star = A^\star$.

\end{definition}

Note that in Definition~\nnn{star}, the equality $(a)^\star = (a)$ could be
replaced by the equality $D = D^\star$ since $(a)^\star = (aD)^\star = aD^\star
= aD = (a)$.

For a brief introduction to star-operations, see \cite[Sections 32 and 34]
{Gilmer}. 
However for a more detailed discussion, see \cite{Jaffard} or \cite{Koch}.

For any star-operation $\star$ and $A, B \in F(D)$, we have that $(A + B)^\star
= (A^\star + B^\star)^\star$ and $(AB)^\star = (AB^\star)^\star = (A^\star
B^\star)^\star$. A star-operation $\star$ {\it distributes over intersections} 
if $(A \cap B)^\star = A^\star \cap B^\star$ for all $A, B \in F(D)$. This is 
easily seen to be equivalent to $(A_1 \cap \cdots \cap A_n)^\star = A_1^\star
\cap \cdots \cap A_n^\star$ for all (integral) $A_1\,,\,.\,.\,.\,, A_n \in
F(D)$. We say that a star-operation $\star$ {\it distributes over arbitrary 
intersections} if $0 \neq \cap A_\alpha \Rightarrow (\cap A_\alpha)^\star =
\cap A_\alpha^\star$ for all $A \in F(D)$. 

A star-operation $\star$ has {\it finite character} if, for each $A \in F(D)$, we have 
that $A^\star = \cup \{ B^\star \mid B \in f(D)$ and $B \subseteq{A} \}$. If $\star$ is a 
star-operation on $F(D)$, we can define a finite character 
star-operation $\star_s$ by $A \longmapsto A^{\star_s} = \cup \{ B^\star \mid 
B \in f(D)$ and $B \subseteq{A} \}$. If $A \in f(D)$, then $A^\star = 
A^{\star_s}$. Three important examples of star-operations are the 
$d$-{\it operation}, the $v$-{\it operation}, and the $t$-{\it operation}.
The $d$-operation is the identity 
operation $A_d = A$ for all $A \in F(D)$. The $v$-operation is defined by 
$A \longmapsto A_v = (A^{-1})^{-1}$, where $A^{-1} = [D:A] = \{ x \in K \mid
xA \subseteq{D} \}$, or equivalently $A \longmapsto A_v = \cap \{ (x) \mid 
A \subseteq{(x)}, 0 \ne x \in K \}$. Throughout this paper, for $A, B \in F(D)$,
$[A:B] = \{ x \in K \mid xB \subseteq{A} \}$, and $(A:B) = \{ x \in D \mid xB
\subseteq{A} \} = [A:B] \cap D$. Note that for any star-operation $\star$, 
$A^\star \subseteq{A_v}$ for all $A \in F(D)$. The $t$-operation is defined by
$A \longmapsto A_t = \cup \{ B_v \mid B \in f(D)$ and $B \subseteq{A} \}$. Note 
that $t = v_s$. Other examples of star-operations will be given and defined in this paper.

Suppose that $\star_1$ and $\star_2$ are two star-operations on $F(D)$, then 
$\star_1 \leq \star_2$ if $A^{\star_1} \subseteq{A^{\star_2}}$ for all $A \in
F(D)$ (or equivalently, $F_{\star_2}(D) \subseteq{F_{\star_1}(D)}$), and 
$\star_1 = \star_2$ if $A^{\star_1} = A^{\star_2}$ for all $A \in F(D)$. The 
star-operations $\star_1$ and $\star_2$ are {\it equivalent} if $A^{\star_1} =
A^{\star_2}$ for all $A \in f(D)$. According to \cite{Gilmer}, for any star-operation 
$\star$, $d \leq \star \leq v$.

A ring $D^\prime$ is said to be an {\it overring} of $D$ if $D
\subseteq{D^\prime} \subseteq{K}$. Let ${\M S}$ be a multiplicatively closed 
set of ideals (or a multiplicative set of ideals) of $D$, i.e, $I, J \in {\M S}$
implies that $IJ \in {\M S}$. For a fractional ideal $A$ of $D$, define 
$A_{\M S} = \{ x \in K \mid xI \subseteq{A}$, for some $I \in {\M S} \}$. Then
$D_{\M S}$ is an overring of $D$, called the ${\M S}$-{\it transform} of $D$,
and $A_{\M S}$ is a fractional ideal of $D_{\M S}$. In \cite{ABrewer}, Arnold 
and Brewer discuss the relationship between the ideal structure of a commutative
ring $R$ and that of the ${\M S}$-transform of $R$. 

Let ${\M S}$ be a multiplicative set of ideals of $D$. The {\it saturation}, 
$\bar {\M S}$, of ${\M S}$ is defined as $\bar {\M S} = \{ A \mid I 
\subseteq{A} \subseteq{D}$, for some $I \in {\M S} \}$. We say that ${\M S}$
is saturated if ${\M S} = \bar {\M S}$. Note that $\bar{\bar{\mathcal{S}}} = \bar 
{\M S}$ so the saturation of a multiplicative set of ideals is saturated. 
It is easily proved that $I, J \in \bar {\M S}$ implies that $IJ \in \bar
{\M S}$. So $\bar {\M S}$ is a multiplicative set of ideals. Clearly 
$A_{\M S} = A_{\bar {\M S}}$. 

\begin{definition}\label{def1.2}

A nonempty set ${\M F}$ of ideals of D is a {\bf localizing system} if 

(1) $I \in {\M F}$ and J an ideal with $I \subseteq{J} \subseteq{D}$ imply 
    that $J \in {\M F}$, and 
    
(2) $I \in {\M F}$ and J an ideal with $(J:i) \in {\M F}$ for all $i \in I$
    imply that $J \in {\M F}$.

\end{definition}  

In Section 2, we investigate when the map $A \longmapsto A_{\M S}$ is a 
star-operation. This question is answered for ${\M S}$ a multiplicative set of
ideals of $D$ and for ${\M S}$ a multiplicative set of finitely generated ideals
of $D$. An obvious necessary condition is that $D_{\M S} = D$. We also determine when
this star-operation has finite character and we give an example where the 
map $A \longmapsto A_{\M S}$ is not a star-operation. 

In Section 3, we characterize when star-operations of the form $A \longmapsto
\cap AD_\alpha$, with $D = \cap D_\alpha$, have finite character.

I give an idea of what I will be working on in the near future in Section 4.


\section{Star-operations and Generalized Quotient Rings}

Recall that a star-operation has finite character if for each $A \in F(D)$, 
$A^\star = \cup \{ B^\star \mid B \in f(D)$ and $B \subseteq{A} \}$.

Note that the $d$-operation and the $t$-operation (defined in the introduction) have 
finite character, but the $v$-operation need not have finite character.

\begin{definition}\label{def1.1}

Let ${\M S}$ be a multiplicatively closed set of ideals of $D$. Then for $A \in F(D)$, 
$A_{\M S} = \{ x \in K \mid xI \subseteq{A}$, for some $I \in {\M S} \}$.

\end{definition}

Recall also that $D_{\M S}$ is an overring of $D$, and if $A$ is a fractional ideal of $D$,
then $A_{\M S}$ is a fractional ideal of $D_{\M S}$.  

\begin{definition}\label{def1.2}

Let D be an integral domain and $\star$ a star-operation on F(D). Then we 
define $\mathcal{S}_{\star} = \{I \in I(D)\mid I^{\star} = D\}$ and 
$\mathcal{S}_{\star}^{f} = \{I \in {\M S}_\star \mid I$ is finitely generated\}.
In particular, $\mathcal{S}_{v} = \{I \in I(D) \mid I_{v} = D$ or 
$I^{-1} = D\}$ and $\mathcal{S}^{f}_{v} = \{I \in \mathcal{S}_{v} \mid {I}$ is
finitely generated \}.

\end{definition}

\begin{definition}\label{def1.3}

Let $I \in F(D)$, then I is said to be a {\bf GV-ideal} (Glaz-Vasconcelos ideal)
if $I_{v} = D$, or equivalently $I^{-1} = D$. A multiplicative set $\mathcal{S}$
is said to be a {\bf GV-system} if $\mathcal{S} \subseteq{\mathcal{S}_{v}}$.

\end{definition}

\begin{proposition}\label{prop1.1}

D = $D_{\mathcal{S}}$ if and only if $\mathcal{S}\subseteq{\mathcal{S}_{v}}$. 

\end{proposition} 

\begin{proof}

$(\ra)$ Let $I \in \mathcal{S}$. Since $I^{-1}I \subseteq{D}$, $I^{-1}
\subseteq{D_{\M S}} = D$. We know that $D \subseteq{I^{-1}}$. Therefore
$I^{-1} = D$. So $I \in \mathcal{S}_{v}$. Hence $\mathcal{S} \subseteq
{\mathcal{S}_v}$.

$(\la)$ Now $D \subseteq{D_{\M S}}$.
Let $x \in D_\mathcal{S}$, then there exists $I \in \mathcal{S} 
\subseteq{{\M S}_v}$ such that $xI \subseteq{D}$. This implies that $x \in 
[ D:I ] = I^{-1} = D$ since $I \in \mathcal{S}_{v}$. Therefore $D_\mathcal{S} 
\subseteq{D}$, and hence $D = D_\mathcal{S}$.
\end{proof}

\begin{proposition}\label{prop1.2}

Let $\mathcal{S}$ be a multiplicative set of nonzero ideals of D. Let $A, B, 
A_1 \,,\,.\,.\,.\,,A_n \in F(D)$ and $a \in K^{\ast}$, then

(1) A $\subseteq{B}$ implies that $A_\mathcal{S} \subseteq{B_\mathcal{S}}$,

(2) A $\subseteq{A_\mathcal{S}}$,

(3) $(aA)_\mathcal{S} = aA_\mathcal{S}$,

(4) $(A_{1} \cap \cdots \cap A_{n})_\mathcal{S} = A_{1_\mathcal{S}} \cap \cdots 
    \cap A_{n_\mathcal{S}}$,

(5) If A and B are integral ideals of D with $A_\mathcal{S} = B_\mathcal{S} = 
    D_\mathcal{S}$, then $(AB)_\mathcal{S} = D_\mathcal{S}$, and
    
(6) The following are equivalent:

\qquad (a) $\mathcal{S} \subseteq{\mathcal{S}_{v}}$,

\qquad (b) $D = D_\mathcal{S}$, and 

\qquad (c) $(a)_\mathcal{S} = (a)$ for all $a \in K^{\ast}$.

\end{proposition}

\begin{proof}

(1) Let $x \in A_{\M S}$, then there exists $I \in {\M S}$ such that $xI 
\subseteq{A} \subseteq{B}$. Therefore $x \in B_{\M S}$ by definition. 
Hence $A_{\M S} \subseteq{B_{\M S}}$.

(2) Let $x \in A$ and $I \in {\M S}$, then $xI \subseteq{AD} \subseteq{A}$.
So $x \in A_{\M S}$, and therefore $A \subseteq{A_\mathcal{S}}$.

(3) Let $x \in (aA)_{\M S}$, then there exists $I \in {\M S}$ such that 
$xI \subseteq{aA}$. Since $K$ is a field and $a \in K^{\ast}$, then 
$\frac{x}{a}I \subseteq{A}$. So $\frac{x}{a} \in A_{\M S}$. Therefore $x \in
aA_{\M S}$. So $(aA)_{\M S} \subseteq aA_{\M S}$.

Suppose $x \in aA_{\M S}$, then $\frac{x}{a} \in A_{\M S}$ since $K$ is a field
and $a \in K^{\ast}$. Therefore there exists $I \in {\M S}$ such that $\frac{x}
{a}I \subseteq{A}$. This implies that $xI \subseteq{aA}$ so $x \in (aA)_{\M S}$.
Therefore $aA_{\M S} \subseteq (aA)_{\M S}$, and hence $(aA)_{\M S} = 
aA_{\M S}$.
 
(4) Since ${A_{1}} \cap \cdots \cap A_{n} \subseteq{A_i}$ for each 
$i = 1,\,.\,.\,.\,,n$, then by (1) $(A_{1} \cap \cdots \cap A_{n})_{\M S} 
\subseteq{A_{i_{\M S}}}$ for each $i = 1,\,.\,.\,.\,,n$. Therefore $(A_{1} \cap \cdots
\cap A_{n})_{\M S} \subseteq {A_{1_\mathcal{S}}} \cap \cdots \cap 
A_{n_\mathcal{S}}$.

Let $x \in A_{1_\mathcal{S}} \cap \cdots \cap A_{n_\mathcal{S}}$, 
then $x \in A_{i_\mathcal{S}}$ for each $i = 1,\,.\,.\,.\,,n$. Therefore there 
exist
$I_{1}$,\,.\,.\,.\,,$I_{n} \in {\M S}$ such that $xI_{1} \subseteq {A_{1}},\,.\,
.\,.\,,xI_{n}
\subseteq {A_{n}}$. So for each $A_{i}$ we have 
$xI_{1} \cdots I_{n} = xI_1 \cdots I_i \cdots I_n \subseteq{xI_i} 
\subseteq{A_i}$. Therefore $xI_1 \cdots I_n \subseteq{A_1 \cap \cdots 
\cap A_n}$. Since ${\M S}$ is multiplicatively closed, $I_{1} \cdots
I_{n} \in {\M S}$. 
Therefore $x \in (A_{1} \cap \cdots \cap A_{n})_
\mathcal{S}$. Thus $A_{1_\mathcal{S}} \cap \cdots \cap A_{n_
\mathcal{S}} \subseteq {(A_{1} \cap \cdots \cap A_{n})_\mathcal{S}}$, and hence
$(A_{1} \cap \cdots \cap A_{n})_\mathcal{S} = A_{1_\mathcal{S}} 
\cap \cdots \cap A_{n_\mathcal{S}}$.

(5) Suppose that $A, B \in I(D)$. Then $AB \subseteq{D}$. Therefore 
$(AB)_{\M S} \subseteq{D_{\M S}}$.

Since $A_{\M S} = D_{\M S}$, $1 \in A_{\M S}$. This implies that there exists
$I \in {\M S}$ such that $1I \subseteq{A}$. Similarly, $B_{\M S} = D_{\M S}$ 
implies that there exists $J \in {\M S}$  such that $1J \subseteq{B}$. So 
$1IJ = 1I1J \subseteq{AB}$. Therefore since $IJ \in {\M S}$, $1 \in 
(AB)_{\M S}$. So $D_{\M S} \subseteq{(AB)_{\M S}}$. Hence $(AB)_{\M S} = 
D_{\M S}$.

(6) $(a) \Leftrightarrow (b)$ This is true by Proposition \ref{prop1.1}.

$(b) \Rightarrow (c)$ Suppose $D = D_{\M S}$, and let $a \in K^{\ast}$.
Then $(a) = aD = aD_{\M S} = (aD)_{\M S} = (a)_{\M S}$. 

$(c) \Rightarrow (b)$ $1 \in K^\ast$. So take $a = 1$. Then by $(c)$, 
$(1)_{\M S} = (1)$. Therefore $D_{\M S} = (1)_{\M S} = (1) = D$.  
\end{proof}

We now state our first theorem which characterizes when the map $A \longmapsto A_{\M S}$,
for $A \in F(D)$ and ${\M S}$ a multiplicative set of finitely generated nonzero ideals, is
a star-operation.

\begin{theorem}\label{thm1.1}

Let $\mathcal{S}$ be a multiplicative set of finitely generated nonzero 
ideals of D. The map $A \longmapsto A_\mathcal{S}$ ($A \in F(D)$) is a 
star-operation if and only if $\mathcal{S} \subseteq{\mathcal{S}_v}$. In this 
case the star-operation $A \longmapsto A_\mathcal{S}$ has finite character. 

\end{theorem} 

\begin{proof}

$(\ra)$ Suppose the map $A \longmapsto A_{\M S}$ is a star-operation, then $D = 
D_{\M S}$. Therefore ${\M S} \subseteq{{\M S}_v}$, by Proposition \ref{prop1.2}. 

$(\la)$ Since $\mathcal{S}$ is a multiplicative set of nonzero ideals, then by
Proposition~\nnn{prop1.2}, $A \subseteq{A_\mathcal{S}}$ (and thus 
$A_\mathcal{S} \subseteq{(A_\mathcal{S})_\mathcal{S}}$), $(xA)_\mathcal{S} = 
xA_\mathcal{S}$, and $A \subseteq{B}$ implies that $A_\mathcal{S} 
\subseteq{B_\mathcal{S}}$. Also by Proposition \ref{prop1.2}, since 
$\mathcal{S} \subseteq{\mathcal{S}_v}$, then $(x)_\mathcal{S} = (x)$ for all
$x \in K^\ast$. So it suffices to show that $(A_\mathcal{S})_\mathcal{S} 
\subseteq{A_\mathcal{S}}$.

Let $x \in (A_\mathcal{S})_\mathcal{S}$, then there exists $I \in \mathcal{S}$
such that $xI \subseteq{A_\mathcal{S}}$. $I$ is finitely generated so suppose
that $I = (x_1,\,.\,.\,.\,,x_n)$. Then for each $x_i \in I, xx_i \in
A_\mathcal{S}$. Therefore for each $x_i$, there exists $A_i \in 
\mathcal{S}$ such that $xx_iA_i \subseteq{A}$. This implies that
$xx_i\prod_{i=1}^{n} A_i \subseteq{xx_iA_i} \subseteq{A}$ for each $x_i$.
So $xI\prod_{i=1}^{n} A_i \subseteq{A}$ for each $i$. Now $I\prod_{i=1}^{n} A_i 
\in 
{\M S}$ since ${\M S}$ is multiplicatively closed. Therefore $x \in A_{\M S}$,
and hence $(A_\mathcal{S})_\mathcal{S} \subseteq{A_\mathcal{S}}$. So the map $A 
\longmapsto A_\mathcal{S}$ is a star-operation. 

We now proceed to show that the map $A \longmapsto A_\mathcal{S}$ has finite 
character (i.e, show that $A_\mathcal{S} = \cup \{ B_\mathcal{S} \mid B \in 
f(D)$ and $B \subseteq{A} \}$).

Let $x \in \cup \{ B_\mathcal{S} \mid B \in f(D)$ and $B \subseteq{A} \}$, then 
$x \in B_\mathcal{S}$ such that $B \in f(D)$ and $B\subseteq{A}$ for some B.
But $B \subseteq{A}$ implies that $B_\mathcal{S} \subseteq{A_\mathcal{S}}$.
Therefore $x \in A_\mathcal{S}$, and hence $\cup \{ B_\mathcal{S} \mid B \in
f(D)$ and $B \subseteq{A} \} \subseteq{A_\mathcal{S}}$.

Let $x \in A_\mathcal{S}$, then there exists $I \in \mathcal{S}$ such that
$xI \subseteq{A}$. Since $I$ is finitely generated, $xI$ is a finitely
generated fractional ideal contained in $A$. So if we let $B = xI$, then
$xI = B$ implies that $x \in B_{\M S}$. So $x \in \cup \{ B_\mathcal{S} \mid 
B \in f(D)$ and $B \subseteq{A} \}$. Therefore $A_{\M S} \subseteq{\cup 
\{ B_\mathcal{S} \mid B \in f(D), B \subseteq{A} \}}$, and hence $A_{\M S}
= \cup \{ B_\mathcal{S} \mid B \in f(D)$ and $B \subseteq{A} \}$. Thus the map
$A \longmapsto A_{\M S}$ has finite character. 
\end{proof} 

Note that the above result is false if the ideals in ${\M S}$ are not assumed to be 
finitely generated as the following example will show. But first, we state a corollary to 
Theorem 2.6

\begin{corollary}\label{cor1.1}

Let ${\M S}$ be a multiplicative set of finitely generated nonzero GV-ideals 
of D. Then the map $A \longmapsto A_{\M S} = A^\star$ is a finite character 
star-operation on D that distributes over finite intersections.

\end{corollary}

\begin{proof}

Since ${\M S}$ is a multiplicative set of GV-ideals, ${\M S} 
\subseteq{{\M S}_v}$ by definition. Also ${\M S}$ consists of finitely 
generated ideals by hypothesis. Therefore by Theorem \ref{thm1.1}, the map
$A \longmapsto A_{\M S} = A^\star$ is a finite character star-operation on $D$.
By Proposition \ref{prop1.2} we have that the map $A \longmapsto A_{\M S}
= A^\star$ distributes over finite intersections.
\end{proof} 

\begin{example}\label{ex2.10a}

Let $D = K[\{ x_n \}_{n = 1}^\infty]$, where $K$ is a field, and let $M = 
(\{ x_n \}_{n = 1}^\infty)$. Then $M$ is a maximal ideal of $D$. Now $M_v = D$.
Let ${\M S} = \{ M^n \}_{n = 1}^\infty$. Now $M_v = D$ implies that $(M^n)_v = D$
for each $n$. So ${\M S} \subseteq{{\M S}_v}$. Therefore $D = D_{\M S}$ by 
Proposition \ref{prop1.2}.
Let $A = (\{ x_n M^n \}_{n = 1}^\infty)$. Now $x_nM ^n \subseteq{A}$ so $x_n
\in A_{\M S}$. Hence $M = (\{ x_n \}_{n = 1}^\infty) \subseteq{A_{\M S}}$. Now
$1 \notin A_{\M S}$ for $1M^n = M^n \nsubseteq{A}$ for any $n$. So since 
$M \subseteq{A_{\M S}} \subsetneq{D}$ and $M$ is a maximal ideal of $D$, 
$M = A_{\M S}$. Therefore $M_{\M S} = (A_{\M S})_{\M S}$. 

Now $M \in {\M S}$, and $1M \subseteq{M}$. Therefore $1 \in M_{\M S}$ and thus
$M_{\M S} = D_{\M S}$. So $(A_{\M S})_{\M S} = M_{\M S} = D_{\M S} = D$.
Hence $A_{\M S} = M \subsetneq{D} = (A_{\M S})_{\M S}$. Therefore in this
example, the map $A \longmapsto A_{\M S}$ is not a star-operation. 

\end{example}

Next, we give a more general characterization of Theorem \ref{thm1.1}. But first, we look 
at some useful results.

Recall that for, ${\M S}$, a multiplicative set of ideals of $D$, the saturation 
of ${\M S}$ is defined by $\bar {\M S} = \{ A \mid I \subseteq{A} \subseteq{D}$,
where $I \in {\M S} \}$. Clearly ${\M S} \subseteq{\bar {\M S}}$ since for any
$I \in {\M S}$, $I \subseteq{I} \subseteq{D}$, which implies that $I \in \bar
{\M S}$.

\begin{definition}\label{def1.4}

A nonempty set ${\M F}$ of ideals of D is a {\bf localizing system} if 

(1) $I \in {\M F}$ and J an ideal with $I \subseteq{J} \subseteq{D}$ imply 
    that $J \in {\M F}$, and 
    
(2) $I \in {\M F}$ and J an ideal with $(J:i) \in {\M F}$ for all $i \in I$
    imply that $J \in {\M F}$.

\end{definition}

\begin{proposition}\label{prop1.3}

Let ${\M S}$ be a multiplicative set of ideals of D. Then for $I \in I(D)$,
$I_{\M S} = D_{\M S}$ if and only if $I \in \bar {\M S}$.

\end{proposition}

\begin {proof}

$(\Rightarrow)$ $I_{\M S} = D_{\M S}$ implies that $1 \in I_{\M S}$. So there
exists $J \in {\M S}$ such that $1J \subseteq{I}$. Therefore $J = 1J 
\subseteq{I} \subseteq{D}$. Thus $I \in \bar {\M S}$.

$(\Leftarrow)$ Suppose that $I \in \bar {\M S}$. Then there exists $J \in 
{\M S}$ such that $J \subseteq{I} \subseteq{D}$. This implies that $1J = J
\subseteq{I}$. Therefore $1 \in I_{\M S}$. So $D_{\M S} \subseteq{I_{\M S}}$.
Hence $I_{\M S} = D_{\M S}$.
\end{proof}

\begin{theorem}\label{thm1.2}

Let ${\M S}$ be a multiplicative set of ideals of D. Then the following are 
equivalent:

(1) $A_{\M S} = (A_{\M S})_{\M S}$ for each $A \in F(D)$,

(2) $I_{\M S} = (I_{\M S})_{\M S}$ for each $I \in I(D)$, and

(3) $\bar {\M S}$ is a localizing system.

\end{theorem}

\begin{proof}

$(1) \Rightarrow (2)$ This is clear since $I(D) \subseteq{F(D)}$.

$(2) \Rightarrow (3)$ Suppose that $I \in \bar {\M S}$ and $I \subseteq{J}
\subseteq{D}$ for some ideal $J$ in $D$. Then by definition, $J \in 
\bar{\bar {\M S}} = \bar {\M S}$. 

Now suppose that $I \in 
\bar {\M S}$ and $J \in I(D)$ with $(J:i) \in \bar {\M S}$ for all $i 
\in I$. Thus $(J : i)_{\M S} = D_{\M S}$ by Proposition \ref{prop1.3}. So 
$(J : i)i \subseteq{J}$ gives $(J : i)_{\M S}(i)_{\M S} \subseteq{J_{\M S}}$
and hence for each $i \in I$, $i \in (i)_{\M S} = (i)_{\M S}D_{\M S}
= (J : i)_{\M S}(i)_{\M S} \subseteq{J_{\M S}}$. This implies that 
$I \subseteq{J_{\M S}}$. Thus $I_{\M S} \subseteq{(J_{\M S})_{\M S}} = 
J_{\M S}$. The last equality follows from the hypothesis. Now by 
Proposition \ref{prop1.3}, $I_{\M S} = D_{\M S}$ so $D_{\M S} = I_{\M S} 
\subseteq{J_{\M S}}
\subseteq{D_{\M S}}$. Hence $J_{\M S} = D_{\M S}$, and so again by 
Proposition \ref{prop1.3}, $J \in \bar {\M S}$. So $\bar {\M S}$ is a localizing system. 

$(3) \Rightarrow (1)$ Suppose that $\bar {\M S}$ is a localizing system and
let $A \in F(D)$. We know that $A_{\M S} \subseteq{(A_{\M S})_{\M S}}$. Let 
$x \in (A_{\M S})_{\M S}$. Then there exists $I \in {\M S}$ such that $xI 
\subseteq{A_{\M S}}$. Let $J = (A : x)$. Then $J \in I(D)$. Now for $i \in I$, 
$xi \in A_{\M S}$. So there exists $J_i \in {\M S}$ such that $xiJ_i 
\subseteq{A}$. So $iJ_i \subseteq{(A : x)} = J$, that is $J_i 
\subseteq{(J : i)}$. Hence $(J : i) \in \bar {\M S}$ for each $i \in I$ since 
$\bar {\M S}$ is a localizing system. So $J \in \bar {\M S}$, again since
$\bar {\M S}$ is a localizing system. Now $xJ \subseteq{A}$ so $x \in 
A_{\bar {\M S}} = A_{\M S}$. Thus $(A_{\M S})_{\M S} \subseteq{A_{\M S}}$. 
Hence $A_{\M S} = (A_{\M S})_{\M S}$.
\end{proof}

We now state the more general characterization of Theorem \ref{thm1.1}

\begin{theorem}\label{thm1.3}

Let ${\M S}$ be a multiplicative set of ideals. Then for $A \in F(D)$, the map
$A \longmapsto A_{\M S}$ is a star-operation if and only if $\bar {\M S}$ is a
localizing GV-system.

\end{theorem}

\begin{proof}

$(\ra)$ Suppose that the map $A \longmapsto A_{\M S}$ is a star-operation for 
$A \in F(D)$. Then $A_{\M S} = (A_{\M S})_{\M S}$. Therefore $\bar {\M S}$ is 
a localizing system by Theorem \ref{thm1.2}.

Let $I \in \bar {\M S}$. Then there exists $J \in {\M S}$ such that $J 
\subseteq{I} \subseteq{D}$. Since the map $ A \longmapsto A_{\M S}$ is a 
star-operation, $(x)_{\M S} = (x)$ for all $x \in K^\ast$. So ${\M S} 
\subseteq{{\M S}_v}$ by Proposition \ref{prop1.2}. 
This implies that $J_v = D$. So $D = J_v \subseteq{I_v} \subseteq{D_v} = D$.
Therefore $I_v = D$ which implies that $I \in {\M S}_v$. Thus $\bar {\M S}
\subseteq{{\M S}_v}$. Hence $\bar {\M S}$ is a localizing GV-system. 

$(\la)$ Since $\bar {\M S}$ is a GV-system, $\bar {\M S} \subseteq{{\M S}_v}$. 
So ${\M S} \subseteq{\bar {\M S}} \subseteq{{\M S}_v}$, and thus $D = D_{\M S}$ 
by Proposition \ref{prop1.2}. 

Since $\bar {\M S}$ is a localizing system, $A_{\M S} = (A_{\M S})_{\M S}$ for
all $A \in F(D)$ by Theorem \ref{thm1.2}. The other conditions for a 
star-operation are satisfied by Proposition \ref{prop1.2}. Hence the map 
$A \longmapsto A_{\M S}$ is a star-operation.   
\end{proof}


\section{Star-operations induced by the localizations of $D$}

In this section, we study the star-operations that are induced by localizations 
$D_{P_\alpha}$ of $D$ at $P_\a$ where $\{P_\a\}$ is a collection of nonzero prime ideals 
of $D$. We then extend our study to the case where the star-operations are induced by 
arbitrary overrings of $D$.

\begin{lemma}\label{lem2.1}

Let $A, B \in F(D)$. Then $(A:B)B \subseteq{A \cap B}$. If $B$ is principal, 
then $(A:B)B = A \cap B$. 

\end{lemma}

\begin{proof}

Let $a \in (A:B)$. Then $aB \subseteq{A}$. Therefore $(A:B)B \subseteq{A}$. 
Suppose $a \in (A:B)$ and $b \in B$. Then $ab \in B$ since $a \in D$. Therefore 
$(A:B)B \subseteq{B}$. Thus $(A:B)B \subseteq{A \cap B}$.  

Suppose $B$ is principal, then $B = (x)$ for some $x \in K^\ast$. Let $a \in
A \cap B$. Then $a \in A$ and $a \in B$. Now $a \in B = (x)$ implies that 
$a = bx$ for some $b \in D$, and $a \in A$ implies that $bx = a \in A$. 
Therefore $b(x) = b(xD) = (bx)D \subseteq{AD} \subseteq{A}$. So $b \in (A:(x))$. 
This implies that $a = bx \in (A:(x))(x) = (A:B)B$, and thus $A \cap B 
\subseteq{(A:B)B}$. Hence if $B$ is principal, then $(A:B)B = A \cap B$. 
\end{proof}

\begin{theorem}\label{thm2.1}

Let A, B $\in F(D)$, $x \in K^{\ast}$ and $\star : F(D) \longrightarrow F(D)$ an 
operation satisfying: 

(1) $A \subseteq{B}$ implies $A^{\star} \subseteq {B^{\star}}$,

(2) $(xA)^{\star} = xA^{\star}$,

(3) $(A \cap B)^{\star} = A^{\star} \cap B^{\star}$,

(4) $A^{\star} = B^{\star} = D$ implies $(AB)^{\star} = D$, and 

(5) $D^{\star} = D$ (or equivalently, $(x)^{\star} = (x)$ for all $x \in 
    K^{\ast}$).
    
Then ${\M S}_{\star} = \{ I \in I(D) \mid I^{\star} = D \}$ is a multiplicative 
set of nonzero ideals of $D$, and $A^{\star} = A_{{\M S}_{\star}}$.    

\end{theorem}

\begin{proof}

Let $I, J \in {\M S}_{\star}$. Then $I^{\star} = D = J^{\star}$. Therefore 
$(IJ)^{\star} = D$ by (4). So $IJ \in {\M S}_\star$, and thus ${\M S}_\star$ is 
a multiplicative set of nonzero ideals of $D$.

Let $x \in A_{{\M S}_{\star}}$. Then there exists $I \in {\M S}_{\star}$ such 
that $xI \subseteq{A}$. Since $I \in {\M S}_\star$, $I^\star = D$. Therefore
$x \in xD = xI^\star = (xI)^\star \subseteq{A^\star}$. The second equality 
follows from (2), and the last containment follows from (1). So 
$A_{{\M S}_{\star}} \subseteq{A^\star}$.

Let $x \in A^\star$, then $(x) \subseteq{A^\star}$. So $(x) = A^\star
\cap (x) = A^\star \cap (x)^\star = (A \cap (x))^\star = [(A:(x))(x)]^\star
= (A:(x))^\star(x)$. The second and third equalities are due to (5) and (3) 
respectively, and the fourth equality is due to Lemma \ref{lem2.1}.
So $(A:(x))^\star = D$, and thus $(A:(x)) \in {\M S}_\star$. Now 
$x(A:(x)) = (x)(A:(x)) = A \cap (x) \subseteq{A}$. So $x \in A_{{\M S}_\star}$, 
and therefore $A^\star \subseteq{A_{{\M S}_\star}}$. Hence $A^\star = 
A_{{\M S}_\star}$.  
\end{proof}   

We now state our first major theorem of this section.

\begin{theorem}\label{thm2.2}

Let $D = \cap D_{P_\alpha}$ where $\{ P_\alpha \}$ is a set of nonzero prime
ideals of D. Let $\star$ be the star-operation defined by $A^\star = 
\cap AD_{P_\alpha}$. Then  the following are equivalent:

(1) The map $A \longmapsto \cap AD_{P_\alpha} = A^\star$ has finite character, 

(2) $A^\star = D$ implies that there exists a finitely generated ideal $B 
\subseteq{A}$ with 

\qquad $B^\star = D$, and

(3) If A is an ideal of D with $A \nsubseteq{P_\alpha}$ for each $\alpha$, then 
there exists a finitely

\qquad generated ideal $B \subseteq{A}$ with 
$B \nsubseteq{P_\alpha}$ for each $\alpha$.

\end{theorem}

\begin{proof}

$(1) \Rightarrow (2)$ Let $A \in I(D)$. Since $\star$ has finite character, 
$A^\star = \cup \{ B^\star \mid B \in f(D)$ and $B \subseteq{A} \}$. So if
$A^\star = D$, then $1 \in D = A^\star = \cup \{ B^\star \mid B \in f(D)$ and 
$B \subseteq{A} \}$. So $1 \in B^\star$ for some $B^\star$ where $B$ is
finitely generated. Therefore $B^\star = D$ for some finitely generated ideal 
$B \subseteq{A}$. (The proof shows that this holds for any finite character 
star-operation.)

$(2) \Rightarrow (1)$ Now $\star$ distributes over finite intersections.
Therefore by Theorem \ref{thm2.1}, for ${\M S}_\star$ a multiplicative subset 
of $D$, $A^\star = A_{{\M S}_\star}$ for all $A \in F(D)$, and ${\M S}_\star 
\subseteq{{\M S}_v}$.

\begin{claim}\label{claim1}

$A_{{\M S}_\star} = A_{{\M S}_\star^f}$ for all $A \in F(D)$.

\end{claim}

We already know that ${\M S}_\star^f 
\subseteq{{\M S}_\star}$. So let $A \in F(D)$. Then $A_{{\M S}_\star^f} 
\subseteq{A_{{\M S}_\star}}$. Let $x \in A_{{\M S}_\star}$. Then $xI 
\subseteq{A}$ for some $I \in I(D)$ such that $I^\star = D$. So by (2), there 
exists a finitely generated ideal $J \subseteq{I}$ with $J^\star = D$. Thus $xJ 
\subseteq{xI} \subseteq{A}$. But $J \in {\M S}_\star^f$. Therefore $x \in 
A_{{\M S}_\star^f}$. This implies that $A_{{\M S}_\star} 
\subseteq{A_{{\M S}_\star^f}}$. Hence $A_{{\M S}_\star} = A_{{\M S}_\star^f}$.

Now $A_{{\M S}_\star} = A_{{\M S}_\star^f}$ for all $A \in F(D)$. 
Therefore for all $A \in F(D)$, $A^\star = A_{{\M S}_\star} =
A_{{\M S}_\star^f}$ by Claim \ref{claim1}, and ${\M S}_\star^f 
\subseteq{{\M S}_\star} \subseteq{{\M S}_v}$. So by Corollary \ref{cor1.1}, the
map $A \longmapsto \cap AD_{P_\alpha} = A^\star$ has finite character.

$(2) \Rightarrow (3)$ Let $A \in I(D)$ such that $A \nsubseteq{P_\alpha}$ for
each $\alpha$. Then $A^\star = D$ by \cite[Theorem 1]{AAStar}. Therefore by (2),
there exists a finitely generated ideal $B \subseteq{A}$ such that
$B^\star = D$. But by \cite[Theorem 1]{AAStar}, $P_\alpha^\star = 
P_\alpha$ for each $P_\alpha \in \{ P_\alpha \}$. So for each $P_\alpha$, 
$P_\alpha^\star = P_\alpha \subsetneq{D} = B^\star$. Suppose that
$B \subseteq{P_\alpha}$ for some $\alpha$, then $B^\star
\subseteq{P_\alpha^\star}$ which is a contradiction. Therefore $B 
\nsubseteq{P_\alpha}$ for each $\alpha$.  

$(3) \Rightarrow (2)$ Suppose that $A \in I(D)$ with $A^\star = D$. Now 
$P_\alpha^\star = P_\alpha \subsetneq{D} = A^\star$ for each $P_\alpha \in 
\{ P_\alpha \}$. The first equality follows from \cite[Theorem 1]{AAStar}.
Suppose that $A \subseteq{P_\alpha}$ for some $\alpha$. Then $A^\star 
\subseteq{P_\alpha^\star}$ which is again a contradiction. So $A 
\nsubseteq{P_\alpha}$ for each $\alpha$. Therefore by (3), there exists a 
finitely generated ideal $B \subseteq{A}$ with $B \nsubseteq{P_\alpha}$ for 
each $\alpha$. But according to \cite[Theorem 1]{AAStar}, this implies that 
$B^\star = D$.
\end{proof}

I am now considering under what conditions Theorem \ref{thm2.2} can be extended to 
star-operations of the form $A \longmapsto \cap AD_\a$ where $D = \cap D_\a$. I have the 
following result.

Let $D_\a$ be an overring of $D$. According to \cite{AAStar}, a star-operation of the form
$A \longmapsto \cap AD_\a$, where $D = \cap D_\a$, is just the special case of the 
star-operation, $\star_\alpha$, defined by the map $A \longmapsto A^\star = 
\cap (AD_\alpha)^{\star_\alpha}$, where each $\star_\alpha$ (the star-operation 
on the corresponding $D_\alpha$) is the identity operation $d$ on $D_\alpha$.
 
\begin {definition}\label{def2.9b} 

A right R-module M is {\bf flat} if and only if for every exact sequence of 
R-modules,
\[ 0 \longrightarrow N^\prime \stackrel{f} \longrightarrow N \stackrel{g}
\longrightarrow N^{\prime\prime} \longrightarrow 0, \] 
the sequence
\[ 0 \longrightarrow M \otimes_{R} N^\prime \stackrel{1_m \otimes f}
\longrightarrow M \otimes_{R} N \stackrel{1_m \otimes g} \longrightarrow 
M \otimes_{R} N^{\prime\prime} \longrightarrow 0 \] 
is exact. 

\end{definition}

\begin{definition}\label{def2.9c} 

An overring T of a ring R is a {\bf flat overring} if T is a flat R-module.

\end{definition}  
   
According to \cite[Theorem 3.3]{LarMc}, a localization is always flat. 

\begin{theorem}\label{thm3.4}

Let $D = \cap D_\alpha$ where each $D_\alpha$ is a flat overring of D. Then
for $A \in F(D)$, the following are equivalent:

(1) The map $A \longmapsto \cap AD_\alpha = A^\star$ has finite character, and

(2) $A^\star = D$ implies that there exists a finitely generated ideal $B 
    \subseteq{A}$ with 
    
    \qquad $B^\star = D$. 
    
\end{theorem}

\begin{proof}

$(1) \Rightarrow (2)$ This follows from Theorem \ref{thm2.2}.

$(2) \Rightarrow (1)$ Suppose that $A^\star = D$ where $A^\star = \cap 
AD_\alpha$. Recall that the operation $A \longmapsto A^\star = \cap AD_\alpha$
is a special case of the star-operation $\star_\alpha$, defined earlier, where 
each $\star_\alpha$ is the identity operation $d$ on $D_\alpha$. So let $A, B
\in F(D_\alpha)$. Then $(A \cap B)^{\star_\alpha} = (A \cap B)_d = A \cap B = 
A_d \cap B_d = A^{\star_\alpha} \cap B^{\star_\alpha}$. Since each $D_\alpha$ 
is flat, $(A \cap B)^\star = A^\star \cap B^\star$ for all $A, B \in F(D)$ by 
\cite[Theorem 2]{AAStar}. That is, $\star$ distributes over intersections. We 
now proceed in the same way as in the proof of Theorem \ref{thm2.2} ($(2) \ra 
(1)$) to conclude that the map $A \longmapsto \cap AD_\alpha = A^\star$ has 
finite character. 
\end{proof}


\section{Future Work}

In his paper, ``Star-Operations Induced by Overrings'', D. D. Anderson characterized when 
a star-operation $\star$ is given by the map $A \longmapsto A^\star = \cap 
AD_{P_\alpha}$, with $D = \cap D_{P_\alpha}$.     
 
I am presently working on characterizing when a star-operation $\star$ is given by the 
map $A \longmapsto A^\star = \cap AD_\alpha$, with $D = \cap D_\alpha$.



\small
\vskip 1pc
{\obeylines
\noindent Sharon \ M. \ Clarke
\noindent Department of Mathematics
\noindent 14 MacLean Hall
\noindent The University of Iowa
\noindent Iowa City, Iowa, 52242
\noindent E-mail: smclarke@math.uiowa.edu
}


\begin{thebibliography}{25}

\bibitem{AAStar} D. D. Anderson, 
{\it Star-operations induced by overrings}, 
Comm. Algebra 16(1988), 2535-2553.
%
\bibitem{AAExam} D. D. Anderson and D. F. Anderson, 
{\it Examples of star-operations on an integral domain}, 
Comm. Algebra 18(1990), 1621-1643.
%
\bibitem{ACook} D. D. Anderson and S. J. Cook, 
{\it Two star-operations and their induced lattices}, 
Comm. Algebra 28(2000), 2461-2475.
%
\bibitem{ABrewer} J. T. Arnold and J. W. Brewer,
{\it On Flat Overrings, Ideal Transforms and Generalized Transforms of a Commutative 
Ring}, Journal of Algebra 18(1971), 254-263.
%
\bibitem{Gilmer} R. Gilmer, 
{\it Miltiplicative Ideal Theory},
Queen's in Pure and Applied Mathematics, vol. 90, Queen's University, Kingston, Ontario, 
Canada, 1992.
%
\bibitem{Koch} F. Halter-Koch,
{\it Ideal {S}ystems, {A}n {I}ntroduction to {M}ultiplicative {I}deal {T}heory},
Monographs and Textbooks in Pure and Applied Mathematics, vol. 211, Marcel Dekker, Inc, 
New York, New York, 1998.
%
\bibitem{Jaffard} P. Jaffard,
{\it Les {S}yst$\grave{e}$mes d'{I}deaux},
Travaux et Recherches Math$\acute{\text{e}}$matiques, vol. IV, Dunod, Paris, 1960.
%
\bibitem{LarMc} M. D. Larsen and P. J. McCarthy,
{\it Multiplicative {T}heory of {I}deals}, 
Pure and Applied Mathematics, vol. 43, Academic Press, Inc, New York and London, 1971.
             

\end{thebibliography}
\end{document}